\newcommand{\rem}[1]{}
\documentclass{amsart}
\usepackage{amsfonts,amssymb,amsmath,amsthm,mathrsfs}
\usepackage{url}
\usepackage[dvips]{epsfig}
\newtheorem{thrm}{Theorem}[section]

\newtheorem{proposition}[thrm]{Proposition}

\newtheorem{remark}[thrm]{Remark}
\theoremstyle{definition}
\newtheorem{definition}[thrm]{Definition}

\begin{document}
\author[C.~A.~Mantica and L.~G.~Molinari]
{Carlo~Alberto~Mantica and Luca~Guido~Molinari}
\address{C.~A.~Mantica: I.I.S. Lagrange, Via L. Modignani 65, 
20161, Milano, Italy -- L.~G.~Molinari (corresponding author): Physics Department,
Universit\`a degli Studi di Milano and I.N.F.N. sez. Milano,
Via Celoria 16, 20133 Milano, Italy.}
\email{carloalberto.mantica@libero.it, luca.molinari@unimi.it}
\title[On conformally recurrent manifolds 
of dimension greater than 4 ]{On conformally  recurrent manifolds 
of dimension greater than 4}
\centerline{Dedicated to the memory of  Dr. Wies{\l}aw Grycak}
{\quad}\\
\begin{abstract} 
Conformally recurrent pseudo-Riemannian manifolds
of dimension $n\ge 5$ are investigated.
The Weyl tensor is represented as a Kulkarni-Nomizu product.
If the square of the Weyl tensor is nonzero, a  covariantly constant symmetric tensor
is constructed, that is quadratic in the Weyl tensor.
Then, by Grycak's theorem, the explicit expression of the traceless part of the Ricci tensor is obtained,
up to a scalar function. The Ricci tensor has at most two distinct
eigenvalues, and the recurrence vector is an eigenvector.\\
Lorentzian conformally recurrent manifolds are then considered. If the square of the Weyl tensor is nonzero,
the manifold is decomposable. A null recurrence vector
makes the Weyl tensor of algebraic type IId or higher in the Bel-Debever-Ortaggio classification, while a
time-like recurrence vector makes the Weyl tensor purely electric.
\end{abstract}
\date{3 march 2016}
\maketitle
\section{Introduction}
Riemannian or pseudo-Riemannian manifolds with a recurrent curvature tensor and generalisations
are the subject of a vast literature.
The recurrent Riemann tensor was first studied in dimension $n=3$ by Ruse in 1949 \cite{Ruse2}
and then by Walker \cite{Walker} (see also Chaki \cite{Chaki2}, Kaigorodov
\cite{Kaigorodov}, Khan \cite{Khan}). Its relationship
with plane waves in general relativity was investigated by Sciama \cite{Sciama}.
Soon after Patterson \cite{Patterson2} introduced Ricci-recurrent spaces.\\
Riemannian manifolds with recurrent Weyl curvature tensor (conformally recurrent manifolds)
were introduced by Adati and Miyazawa \cite{Adati}, and generalised to pseudo-Riemannian manifolds
by Derdzi\'nski \cite{Derdzinski}, Roter \cite{Roter00,Roter1,Roter2,Roter3} and others, for instance by Suh and Kwon \cite{SuhKwon}.
Mc Lenaghan and Leroy \cite{LenaghLer} and Mc Lenaghan
and Thompson \cite{LenaghThomps} considered Lorentzian manifolds (space-times)
with complex recurrent Weyl tensor. They showed that such
spaces belong to Petrov types D or N, and obtained the
expression of the metric in the case of real recurrence vector. Conformally
recurrent space-times were also studied by Hall \cite{Hall74,Hall77}.
\begin{definition}
A pseudo-Riemannian manifold $(\mathcal M,g)$ of dimension $n\ge 4$ is said to be conformally recurrent
if its Weyl conformal curvature tensor\footnote{The components
of the Weyl tensor are \cite{Postnikov}:
$$ C_{jkl}{}^m = R_{jkl}{}^m +\frac{1}{n-2}(\delta_j{}^mR_{kl} -\delta_k{}^m
R_{jl}+R_j{}^mg_{kl} - R_k{}^m g_{jl})\\
 - R \frac{\delta_j{}^m g_{kl}-\delta_k{}^m g_{jl}}{(n-1)(n-2)} $$
where $R_{kl} =-R_{mkl}{}^m$ is the Ricci tensor and $R=g^{ij}R_{ij}$ is
the curvature scalar.}
satisfies the condition:
\begin{align}
\nabla_i C_{jklm} = \alpha_i C_{jklm} \label{eq1.1}
\end{align}
with a non-zero covector field $\alpha_i$ named recurrence vector.
\end{definition}
The definition has two straightforward well known consequences. The first
one is an equation for the recurrence vector:
$\nabla_i C^2 = 2 \alpha_i C^2$, where $C^2=C_{jklm}C^{jklm}$.
Then, if $C^2\neq 0$ at a point of $\mathcal M$,  on a coordinate domain $U$ of this point we have
\begin{align}
\alpha_i = \tfrac{1}{2}\nabla_i \log |C^2|. \label{alphaC2}
\end{align}
Another consequence is the identity:
\begin{align}
[\nabla_i, \nabla_j]\,C_{klmn}= (\nabla_i\alpha_j-\nabla_j\alpha_i)C_{klmn}.
\label{eq3}
\end{align}
If the 1-form $\alpha_i$ is closed ($\nabla_i\alpha_j-\nabla_j\alpha_i=0$) then $[\nabla_i, \nabla_j]\,C_{klmn}= 0$, i.e.
the manifold is Weyl semi-symmetric \cite{Des90,Des91,JHPTV,DGHSaw}. 

We recall that a manifold is semi-symmetric if $[\nabla_i, \nabla_j]R_{klmn}=0$ (Cartan \cite{Cartan}, Szab\'o \cite{Szabo}).
Derdzi\'nski and Roter proved that every non-conformally flat, non-locally symmetric semi-Riemannian manifold of dimension $n\ge 4$ 
with parallel Weyl tensor is semi-symmetric \cite{Der78} (theorem 9).
The Riemannian semi-symmetric manifolds were classified by Szab\'o \cite{Szabo}.\\
Semi-symmetry implies Weyl semi-symmetry: in fact from $[\nabla_i,\nabla_j]R_{klmn}=0$ we get $[\nabla_i,\nabla_j]R_{lm}=0$ and then
$[\nabla_i,\nabla_j]C_{klmn}=0$.
In general the converse is not true: a counter example in dimension $n=4$ was given by Derdzi\'nski \cite{Der81} (see also \cite{Des90}).
Equivalence between semi-symmetry and Weyl semi-symmetry holds on non-conformally flat $n=4$ warped product manifolds 
(Deszcz \cite{Des91}, theorem 3 and corollary 1). It was established by Grycak for non-conformally flat pseudo-Riemannian manifolds 
of dimension $n\ge 5$ \cite{Grycak} (see also \cite{DG}). In $n=4$  it holds for Lorentzian manifolds (Eriksson and Senovilla \cite{Eriksson}).

In $n\ge 5$, by the aforementioned equivalence, Weyl semi-symmetry implies Ricci
semi-symmetry, 
\begin{align}
[\nabla_i,\nabla_j]R_{kl}=0 \label{4a}
\end{align}
It seems that \eqref{4a} originally appeared in RoterÕs paper \cite{WR1} (Lemma 3, eq.  10).
A summation
over cyclic permutations of indices $ijk$ in \eqref{4a} gives the algebraic property (Defever and Deszcz \cite{Defever}):
\begin{align}
 R_{im} R_{jkl}{}^m +R_{jm} R_{kil}{}^m +R_{km} R_{ijl}{}^m =0. \label{eq2.1RR}
\end{align}
This equation originally appeared in Roter's paper on conformally symmetric spaces (\cite{Rot72} lemma 1). In \cite{ManMol2012} the 
property was generalized by us, to define the notion of K-compatible tensor, where K is Riemann's or Weyl's or a 
generalized curvature tensor (then Eq.~(\ref{eq2.1RR}) states that the Ricci tensor is Riemann compatible). Geometric implications
of compatibility were explored for the Riemann tensor \cite{RCT} and the
Weyl tensor \cite{WCT}. In particular we showed that Weyl and Riemann compatibility are equivalent for
the Ricci tensor.
Manifolds whose Ricci tensor is Weyl compatible are
termed ``Weyl compatible manifolds'' \cite{ManSuh2014}.
Several such manifolds, which include the Robertson-Walker space-time, were discussed by Deszcz et al. \cite{Deszcz2013}; another example is G\"odel's metric (\cite{Deszcz2014}, th.2). \\

In this paper we present new results on conformally recurrent manifolds with dimension $n\ge 5$.
In Section 2 we specialize the second Bianchi identity for the Weyl tensor and an identity by Lovelock to
algebraic ones, and show that in $n\ge 5 $ the tensor $\alpha_i\alpha_j$ is Weyl compatible; also the Ricci tensor is such,
when $\alpha $ is closed. The Weyl tensor is represented as a Kulkarni-Nomizu product in terms of
the recurrence vector $\alpha_i$ and the symmetric tensor $E_{il}= C_{ijkl}\alpha^j\alpha^k/\alpha^2 $ (where $\alpha^2=\alpha^k\alpha_k$).\\
In Section 3 we introduce the tensor $h_{ij}=C_i{}^{klm}C_{jklm}/C^2$.
Among other properties, we show that it is covariantly constant, $\nabla_k h_{ij}=0$.  Therefore, this tensor
gives an explicit realization of Grycak's theorem \cite{Grycak1981}, namely that the traceless part of the Ricci tensor is
proportional to the traceless part of $h_{ij}$ via a scalar function, Eq.~(\ref{Grycak}).
This means that, up to a scalar function, the traceless part of the Ricci tensor is determined
by the Weyl tensor. Proportionality implies unicity of $h_{ij}$, that the Ricci tensor
cannot have more than 2 distinct eigenvalues, and that the recurrence vector is an eigenvector.\\
In Sections 4 and 5 we focus on Lorentzian conformally recurrent manifolds, $n\ge 5$:
the existence of $h_{ij}$ with zero covariant derivative implies that a Lorentzian conformally recurrent manifold is
either decomposable or it admits a null covariantly constant vector field \cite{Senovilla, Hall91}.
If $C^2\neq 0$ the manifold is
decomposable. A Brinkmann or pp-wave metric requires $C^2=0$, with the characterization
given by Galaev \cite{Galaev}. \\
With null recurrence vector $\alpha^i\alpha_i=0$,
the algebraic type of the Weyl tensor is at least II$_d$ in the high-dimensional Bel-Debever-Ortaggio
classification \cite{Ortaggio09}. This extends the result of
Mc Lenaghan and Leroy \cite{LenaghLer} valid in $n=4$.
If $\alpha^2 < 0$, Weyl's tensor is purely electric, according to the definition given in \cite{HeOrWy2013}.\\
In Section 6 we construct a simple example of conformally recurrent manifold.

We assume that the manifolds are smooth, connected, Hausdorff,
with non-degenerate metric ($n$-dimensional pseudo-Riemannian manifolds)
and that  $\nabla_j g_{kl} = 0$. Where necessary we specialize to a
metric with signature $n-2$ (Lorentzian manifolds).

\section{\bf A representation of the Weyl tensor}
Our presentation largely relies on two general identities for the Weyl tensor
on a pseudo-Riemannian manifold. The first one is eq.(3.7) in \cite{Adati} (see also \cite{ManSuh2013}):
\begin{align}
\nabla_i C_{jklm} + \nabla_j C_{kilm} + \nabla_k C_{ijlm}
= \frac{1}{n-3} \nabla_p
(g_{jm}  C_{kil}{}^p + g_{km} C_{ijl}{}^p \nonumber \\
+ g_{im}  C_{jkl}{}^p
+ g_{kl}  C_{jim}{}^p + g_{il} C_{kjm}{}^p + g_{jl}  C_{ikm}{}^p ).
\end{align}
\noindent
On a conformally recurrent manifold with recurrence vector $\alpha_i$, the identity becomes algebraic:
\begin{align}
\alpha_i C_{jklm} +\alpha_j C_{kilm} +\alpha_k C_{ijlm}
= \frac{\alpha^p}{n-3}  ( g_{jm} C_{kilp} + g_{km} C_{ijlp} \nonumber \\
+ g_{im} C_{jklp} + g_{kl} C_{jimp} + g_{il} C_{kjmp} + g_{jl} C_{ikmp} ) \label{eq2.13}
\end{align}
\begin{remark}{\quad}\\
1) Adati and Miyazawa \cite{Adati} considered the additional condition
$\alpha^m C_{jklm}=0$. Then $\alpha_i C_{jklm} +\alpha_j C_{kilm} +\alpha_k C_{ijlm} =0$ and
contraction with $\alpha^i$ gives $\alpha^2 C_{jklm}=0$, i.e. either $\alpha^2=0$ or $C_{jklm}=0$.\\
2) Equation~(\ref{eq2.13}) holds on manifolds more general than conformally recurrent, with the property
$\nabla_i C_{jklm} + \nabla_j C_{kilm} + \nabla_k C_{ijlm} = \alpha_i C_{jklm} + \alpha_j C_{kilm} + \alpha_k C_{ijlm} $.
It implies $\nabla_m C_{jkl}{}^m =\alpha_m C_{jkl}{}^m$. They are studied in \cite{ManSuh2014b}.
%\cite{ManSuh2012,ManSuh2013,ManSuh2014b}.
\end{remark}
The second identity results from two identities by Lovelock, relating the Riemann, the Ricci and the Weyl tensors
on a pseudo-Riemannian manifold \cite{ManMol2011,RCT}:
\begin{align}
& \nabla_i \nabla_m C_{jkl}{}^m + \nabla_j \nabla_m C_{kil}{}^m +
\nabla_k \nabla_m C_{ijl}{}^m   \nonumber \\
& = - \frac{n-3}{n-2}  ( R_{im}  C_{jkl}{}^m + R_{jm} C_{kil}{}^m +
R_{km} C_{ijl}{}^m ) \label{eq2.6} 
\end{align}

\begin{proposition} \label{theorem 3.3}
Let $\mathcal M$, $n\ge 5$, be a pseudo-Riemannian conformally recurrent manifold with recurrence
vector $\alpha_i$.  The tensors $\alpha_i\alpha_j$, $\nabla_i\alpha_j +\nabla_j\alpha_i + R_{ij}$
and $\nabla_i\alpha_m -\nabla_m\alpha_i +\frac{n-4}{n-2}\, R_{im}$ are Weyl compatible:
\begin{align}
 \alpha_i\alpha_m C_{jkl}{}^m +\alpha_j\alpha_m C_{kil}{}^m +\alpha_k \alpha_mC_{ijl}{}^m = 0
\label{eq2.17}\\
\left[ \nabla_i\alpha_m +\nabla_m\alpha_i + R_{im} \right] C_{jkl}{}^m +\text{cyclic} =0\\
\left[ \nabla_i\alpha_m -\nabla_m\alpha_i +\frac{n-4}{n-2}\, R_{im}\right]C_{jkl}{}^m+\text{cyclic} =0
\end{align}
{\rm cyclic} stands for sum of cyclic permutations on $ijk$.
Moreover, if the recurrence vector is closed, then
$\nabla_i\alpha_j$ and  $R_{ij}$ are Weyl compatible, and
the Ricci tensor is also Riemann compatible, Eq.~(\ref{eq2.1RR}).
\begin{proof}
Contraction of Eq.~(\ref{eq2.13}) with $\alpha^m$ cancels some terms and, with some algebra, results in
Weyl compatibility for $\alpha_i\alpha_j$. By the recurrence property it is:
$\nabla_i \nabla_m C_{jkl}{}^m = \nabla_i(\alpha_m C_{jkl}{}^m)
=(\nabla_i\alpha_m +\alpha_i\alpha_m ) C_{jkl}{}^m $. Equation~(\ref{eq2.6}) becomes:
\begin{align}
\left[ \nabla_i\alpha_m +\alpha_i\alpha_m +\frac{n-3}{n-2}R_{im} \right]C_{jkl}{}^m +
\text{cyclic} =0.\nonumber
\end{align}
The second term is zero by Eq.~(\ref{eq2.17}).
The covariant divergence $\nabla^m$ of Eq.~(\ref{eq2.13}), and (\ref{eq2.17}) give:
$$ \left [ \nabla^m \alpha_i - \frac{1}{n-3}\nabla_i\alpha^m \right ] C_{jklm} + \text{cyclic} = \frac{(\nabla^m\alpha^p)}{n-3} (g_{il}C_{jkmp}-g_{jl}
C_{kimp}-g_{kl}C_{ijmp})$$
Contraction with $g^{kl}$ gives: $ \left [ \nabla_m \alpha_k - \frac{1}{n-3}\nabla_k\alpha_m \right ] C_{ij}{}^{km} = -\frac{n}{n-3}(\nabla^m\alpha^p) C_{ijmp}$ i.e. $(\nabla_m \alpha_k)C_{ij}{}^{km} =0$. We then remain with:
\begin{align}
 \left [ \nabla_m \alpha_i - \frac{1}{n-3}\nabla_i\alpha_m \right ] C_{jkl}{}^m + \text{cyclic} = 0.\nonumber
 \end{align}
Linear combinations of the two cyclic identities produce the final statements.
\end{proof}
\end{proposition}
For $n\ge 5$  and $\alpha^2\neq 0$, a representation of the Weyl tensor can be derived in terms of the recurrence vector
and of the traceless symmetric tensor
\begin{align}
E_{il}= \frac{\alpha^j\alpha^k}{\alpha^2} C_{ijkl}. \label{eq2.15a}
\end{align}
Note that $E_{il}\alpha^l=0$.
%The covariant derivative is evaluated:
%\begin{align}
% \nabla_r E_{ij} =\alpha_r E_{ij} - (\alpha_i E_{jk}+\alpha_j E_{ik}) \frac{\nabla_r\alpha^k}{\alpha^2}. \label{E2}
%\end{align}
%
\begin{thrm}\label{Th3.6}
On a conformally recurrent manifold of dimension $n\ge 5$ with $\alpha^2\neq 0$,
the Weyl tensor has the form:
\begin{align}
 C_{jklm} = \frac{1}{n-3}
 \left [ g_{mk} E_{jl} - g_{mj} E_{kl} + g_{jl}
E_{km}  - g_{kl} E_{jm} \right ] \nonumber  \\
-\frac{n-2}{n-3} \left [ \frac{\alpha_k\alpha_m}{\alpha^2} E_{jl} - \frac{\alpha_j \alpha_m}{\alpha^2}  E_{kl}  +
\frac{\alpha_j \alpha_l}{\alpha^2} E_{km} - \frac{\alpha_k\alpha_l}{\alpha^2} E_{jm}\right ]. \label{eq2.20}
\end{align}
\begin{proof}
Contraction of Eq.~(\ref{eq2.13}) with $\alpha^i$ gives
 \begin{align}
\alpha^2 C_{jklm} +\alpha^p(\alpha_j C_{lmkp} +\alpha_k C_{mljp})
- \frac{\alpha^p}{n-3}  (\alpha_m C_{jklp} + \alpha_l C_{kjmp}) \nonumber \\
=\frac{\alpha^p\alpha^i}{n-3}  (- g_{jm} C_{iklp} + g_{km} C_{ijlp}
- g_{kl} C_{ijmp} + g_{jl} C_{ikmp} ),\nonumber
\end{align}
symmetrization in the exchange of the pairs $jk$ and $lm$ gives
\begin{align}
\alpha^2 C_{jklm} = -\frac{1}{2}\frac{n-2}{n-3}\,\alpha^p\,(\alpha_j C_{lmkp} +\alpha_k C_{mljp}
+ \alpha_m C_{kjlp} + \alpha_l C_{jkmp})  \nonumber\\
 +\frac{\alpha^p\alpha^i}{n-3}  ( -g_{jm} C_{iklp} + g_{km} C_{ijlp}
- g_{kl} C_{ijmp} + g_{jl} C_{ikmp} ), \nonumber \\
= - \frac{1}{2}\frac{n-2}{n-3}\alpha^p (\alpha_j C_{lmkp} +\alpha_k C_{mljp} +
\alpha_m C_{kjlp} + \alpha_l C_{jkmp})\nonumber \\
+\frac{\alpha^2}{n-3}  ( -g_{jm} E_{kl}
+ g_{km} E_{jl} - g_{kl} E_{jm} + g_{jl} E_{km} ).\nonumber
\end{align}
Contraction of Eq.~(\ref{eq2.17}) with $\alpha^i$ gives
$ \alpha^p C_{jklp} = \alpha_j E_{kl}- \alpha_k E_{jl}  $, and Eq.~(\ref{eq2.20}) is obtained.
\end{proof}
\end{thrm}
\begin{remark} The representation Eq.~(\ref{eq2.20}) of the Weyl tensor is the Kulkarni-Nomizu product
\cite{Besse} of the tensors $\frac{1}{n-3}[g_{jk}-(n-2)
\alpha_j\alpha_k/\alpha^2]$ and $E_{ml}$.
\end{remark}
\noindent
Contraction of Eq.~(\ref{eq2.20}) by $C^{jklm}$ gives $C^2 = 4 \,\frac{n-2}{n-3}\, E_{ij}E^{ij}$.
% \label{Esquare}\end{align}
%
\begin{proposition}
On a conformally recurrent manifold $n\ge 5$ with $\alpha^2\neq 0$, the tensor $E_{ij}$ is Weyl compatible:
\begin{align}
 E_{im}C_{jkl}{}^m + E_{jm}C_{kil}{}^m + E_{km} C_{ijl}{}^m =0.
 \end{align}
\begin{proof}
From  the representation Eq.~(\ref{eq2.20}) we evaluate
\begin{align}
 E_{im}C_{jkl}{}^m = \frac{1}{n-3}(E_{ik}E_{jl} -E_{ij}E_{kl}+g_{jl}E_{km}E_i{}^m - g_{kl}
E_{im}E_j{}^m)\nonumber\\
+\frac{n-2}{n-3} \left[ \frac{\alpha_k\alpha_l}{\alpha^2} E_{im}E_j{}^m - \frac{\alpha_j\alpha_l}{\alpha^2} E_{km}E_i{}^m\right];\nonumber
\end{align}
the sum over cyclic permutations of indices $ijk$ cancels all terms in the right-hand side, and
Weyl compatibility is proven.
%Contraction of indices $i$ and $j$ gives:
%$$E_{jm}C^j{}_{kl}{}^m = \frac{1}{n-3} \left[ 2E_{km}E_l{}^m g_{kl}E^2+ (n-2) \frac{\alpha_k\alpha_l}{\alpha^2}
%E^2\right] $$
\end{proof}
\end{proposition}
\section{\bf A realization of a theorem by Grycak, and the traceless part of the Ricci tensor}
In this section we study conformally recurrent manifolds  with $C^2\neq 0$ and dimension $n\ge 5$. By Eq.~(\ref{alphaC2}) the recurrence vector is closed, $\nabla_i\alpha_j=\nabla_j\alpha_i$.
We define the symmetric tensor
\begin{align}
h_{ij} = \frac{1}{C^2} C_i{}^{klm}C_{jklm}. \label{tensorh}
\end{align}
and assume that it is not proportional to the metric tensor (in $n=4$ the tensor is trivial: $h_{ij} = \tfrac{1}{4} g_{ij}$).
The trace is $g^{ik}h_{ik} = h^k{}_k=1$. Recurrence of the Weyl tensor
implies that the tensor is covariantly constant:
\begin{align}
  \nabla_j h_{kl}= 0.
 \end{align}
\begin{proposition}
Let $\mathcal M$, $n>4$, be a conformally recurrent pseudo-Riemannian manifold with non-zero scalar $C^2$.
The tensor $h_{ij}$ defined by \eqref{tensorh} has the following properties:\\
1) the recurrence vector is an eigenvector of $h$:
\begin{align}
h^i{}_j  \alpha^j=  \frac{(n-3)}{2(n-2)}\, \alpha^i \label{eq*3}
\end{align}
2) the tensor $h_{ij} $ is Riemann compatible:
\begin{align}
h_{im}R_{jkl}{}^m +h_{jm}R_{kil}{}^m + h_{km}R_{ijl}{}^m  =0,
 \end{align}
it commutes with the Ricci tensor, and it is Weyl compatible.
\begin{proof}
Contraction of Eq.~(\ref{eq2.13}) with $C^{qklm}$ gives the identity:
\begin{align}
 \alpha_i h_j{}^q -\alpha_j h_i{}^q = -\frac{\alpha^p}{C^2} \left[ C_{jilm}C_p{}^{qlm} + \frac{2}{n-3}
\left( C_{iklp} C^{qkl}{}_j - C_{jklp}C^{qkl}{}_i\right)\right]\label{equat2}
\end{align}
Contraction with $\delta_q{}^i$ gives the first result.
The property $[\nabla_i,\nabla_j] h_{kl}=0$ is, by the Ricci identity,
$R_{ijk}{}^m h_{ml}+R_{ijl}{}^m h_{km} = 0$.
Summation on cyclic permutations of $ijk$ and the first Bianchi identity give Riemann compatibility. This implies Weyl compatibility \cite{WCT}.
Contraction with $g^{jl}$ gives
$h_{im} R_k{}^m -R_i{}^m h_{km} =0$.
\end{proof}
\end{proposition}
The explicit expression \eqref{tensorh} of a nontrivial covariantly constant tensor $h_{ij}$ gives,  in $n>4$, 
the realisation of an interesting theorem by Grycak, that is here
recalled:
\begin{thrm} (Grycak, \cite{Grycak1981} theorem 1)
Let $\mathcal M$ be a conformally recurrent  manifold $n\ge 4$ that is neither conformally flat nor recurrent,
whose recurrence vector $\alpha^i$ is locally a gradient. If $\mathcal M$ admits a symmetric parallel
tensor $h_{ij}$ that is not multiple of the metric, then $
 \left( R_{ij} - \frac{1}{n} R g_{ij} \right)=G \left (h_{ij}-\frac{1}{n} h^k{}_k g_{ij}\right ) $,
$G$ being a scalar function.
\end{thrm}
Here the statement is made explicit:
\begin{align}
\left( R_{ij} - \frac{R}{n}  g_{ij} \right)=G \left (\frac{1}{C^2} C_{iklm}C_j{}^{klm}-\frac{1}{n} g_{ij}\right )
 \label{Grycak}
 \end{align}

\begin{remark}
A covariant derivative shows that the traceless part of the Ricci tensor is recurrent:
\begin{align}
\nabla_k \left (R_{ij} - \frac{R}{n}g_{ij}\right ) = \frac{\nabla_k G}{G}\left(R_{ij} - \frac{R}{n}g_{ij}\right )
\end{align}
\end{remark}
The tensor $h^2_{ij}=h_{ik}h^k{}_j$ is symmetric and has zero covariant derivative. 
Grycak's theorem implies proportionality of the traceless parts of  $h^2_{ij}$ and $h_{ij}$:
\begin{align}
h^2_{ij} - \frac{(h^2)^k{}_k}{n}g_{ij} = H \left (h_{ij}-\frac{1}{n} g_{ij}\right ) \label{hsqh}
\end{align}
where $(h^2)^k{}_k = g^{jk} h^2{}_{jk}$ and $H$ is a scalar (a covariant derivative of both sides
of \eqref{hsqh} gives $\nabla_i H=0$).
\begin{thrm} Let $\mathcal M$, $n\ge 4$, be a conformally recurrent pseudo-Riemannian 
manifold with non-zero scalar $C^2$.
The Ricci tensor has at most two distinct eigenvalues:
\begin{align}
\begin{cases}
\mu=\frac{R}{n}+ G \frac{(n-1)(n-4)}{2n(n-2)} & \text{with multiplicity $n_h$}\\
\mu' =\frac{R}{n}-G\frac{(n-1)(n-4)n_h}{2n(n-2)(n-n_h)} & \text{with multiplicity $n-n_h$.}
\end{cases}
\end{align}
and the spectral decomposition $R^i{}_j = \mu P^i{}_j + \mu' (\delta^i{}_j-P^i{}_j)$ with
projector
\begin{align}
P^i{}_j = \frac{2(n-2)(n-n_h)}{(n-1)(n-4)} h^i_j +  \frac{n_h(n-3)-2(n-2)}{(n-1)(n-4)} \delta^i{}_j \label{pij}
\end{align}
on a submanifold of dimension $n_h$.
The curvature scalar is $R=n_h\mu +(n-n_h)\mu' $.
\begin{proof}
Equation~(\ref{hsqh}) implies that $h_{ij}$ has at most two eigenvalues and that they solve the equation $\lambda^2 -H\lambda -\frac{1}{n} [(h^2)^k_k-H]=0$. They are:
$h=\frac{n-3}{2(n-2)}$ and $h'=H-h$. Let the corresponding eigenspaces have dimensions $n_h$ and $n-n_h$. The first eigenspace contains the recurrence vector. It is:
\begin{align}
  & n_h h + (n-n_h) h' =1\nonumber\\
  & n_h h^2 + (n-n_h) h^{\prime 2} = (h^2)^k{}_k\nonumber
 \end{align}
Let $P^i{}_j$ be the projector on the eigenspace with eigenvalue $h$. Then $P^k{}_k =n_h$, and
$h^i{}_j= h P^i{}_j + h' (\delta^i{}_j -P^i{}_j)$. From this relation the expression Eq.~(\ref{pij}) is obtained.\\
The algebraic constraints are inherited by the Ricci tensor via Grycak's relation. The Ricci tensor
has two eigenvalues $\mu$ and $\mu'$:
\begin{align}
\begin{cases}
\mu=\frac{R}{n}+ G \left( h -\frac{1}{n}\right) & \text{with multiplicity $n_h$}\\
\mu' =\frac{R}{n}+G\left( h' -\frac{1}{n}\right)& \text{with multiplicity $n-n_h$.}
\end{cases}
\end{align}
which are evaluated. The trace of the Ricci tensor is $n_h\mu+(n-n_h)\mu'$.
\end{proof}
\end{thrm}

\begin{remark}{\quad}\\
1) The case where the second eigenvalue of $h_{ij}$ is zero ($h'=0$) implies $n_h=2+\frac{2}{n-3} =$
integer, therefore it may only occur in $n=5$, with $n_h=3$.\\
2) The case where the eigenvalue $\mu $ of the Ricci tensor is non-degenerate ($n_h=1$)
gives, with simple algebra,
\begin{align}
 R_{ij}-\frac{R}{n}g_{ij}= \frac{G}{2} \frac{n-4}{n-2}\left[ \frac{\alpha_i\alpha_j}{\alpha^2} - \frac{g_{ij}}{n} \right].
 \end{align}
Since $h_i{}^j\alpha_j=h\alpha_i$ implies $h_i{}^j\nabla_r\alpha^i =h\nabla_r\alpha_j$, the recurrence vector
is recurrent: $\nabla_r \alpha_j= K \alpha_r\alpha_j$ where $K$ is a number. The manifold is quasi-Einstein  i.e.
there is a scalar field $\lambda $ such that ${\rm rank} (R_{ij}-\lambda g_{ij})\le 1$ (see \cite{Deszcz2014,Chaki2000,Chojnacka,DGJZ}).\\
3) In $n>4$ the eigenvalues $h$ and $h'$ are different. Patterson proved that 
a pseudo-Riemannian manifold is a product manifold if there exists a symmetric recurrent tensor that has at least two different eigenvalues 
at any point \cite{Patterson51} (see also \cite{Duggal} page 96). In the present case, the recurrent tensor is $C^2 h_{ij}$. It follows that
a conformally recurrent manifold with $C^2\neq 0$ and $n>4$ is a product manifold.
\end{remark}

\begin{remark}{\quad}\\
From the Ricci identity we have $R_{jklm}\alpha^m =[\nabla_j,\nabla_k]\alpha_l$ and 
contraction  with $g^{jl}$ gives:
$R_{km} \alpha^m = [\nabla_l,\nabla_k]\alpha^l = \mu \alpha_k$. If $\nabla_i\alpha_j=\nabla_j\alpha_i$ this corresponds to the equation:
\begin{align}
\nabla^2 \alpha_k - \nabla_k \nabla_m\alpha^m  = \mu \,\alpha_k
\end{align}
\end{remark}

The divergence $\nabla^i$ of \eqref{Grycak} gives a linear relation among gradients:
$\frac{n-2}{2n} \nabla_j R = (h^m{}_j - \frac{1}{n} \delta^m{}_j)\nabla_m G$. With some effort we get rid of the tensor $h^m{}_j$:

\begin{proposition}
Let $\mathcal M$, $n>4$, be a conformally recurrent pseudo-Riemannian manifold with
non-zero scalar $C^2$. The function G, defined by \eqref{Grycak}, satisfies the following
equation:
\begin{align}
\nabla_j G = \frac{n(n-2)^2}{(n-1)^2(n-4)} \left [  \frac{\alpha_j \alpha_m}{\alpha^2}-\frac{g_{jm}}{n}    \right ]\nabla^m R
\end{align}
\begin{proof}
Consider the expression for the covariant divergence of the Weyl tensor:
\begin{align}
\nabla^m C_{jklm} = -\frac{n-3}{n-2}\left[ \nabla_jR_{kl}-\nabla_kR_{jl} -\frac{1}{2(n-1)} (g_{kl}\nabla_j R -g_{jl}\nabla_k R)\right] \label{divC}
\end{align}
Contraction with $\alpha^l$ gives
\begin{align}
0=\alpha^l (\nabla_jR_{kl}-\nabla_kR_{jl} ) -\tfrac{1}{2(n-1)} (\alpha_k\nabla_j R -\alpha_j\nabla_k R)\label{80}
\end{align}
By Eq.~(\ref{Grycak}) it is
\begin{align}
\nabla_j R_{kl}= \frac{1}{n}g_{kl} \nabla_j R + (\nabla_j G) \left(h_{kl}-\frac{1}{n}g_{kl}\right) \label{82}
\end{align}
Then:
$$\alpha^l \nabla_j R_{kl} - \frac{\alpha_k \nabla_j R}{2(n-1)}
= \frac{\alpha^k}{2n(n-1)(n-2)} \left[
(n-2)^2 \nabla_j R  + (n-1)^2(n-4) \nabla_j G\right ]$$
The expression is inserted in Eq.~(\ref{80}):
\begin{align}
0=(n-2)^2 (\alpha_k \nabla_j R - \alpha_j \nabla_k R) + (n-1)^2(n-4)(\alpha_k \nabla_j G -
 \alpha_j \nabla_k G)\nonumber
\end{align}
Contract with $\alpha^k$:
\begin{align}
0=(n-2)^2 (\alpha^2 \nabla_j R - \alpha_j \alpha^k \nabla_k R) + (n-1)^2(n-4)(\alpha^2 \nabla_j G -
 \alpha_j \alpha^k\nabla_k G)\nonumber
\end{align}
A simplification occurs: contraction of Eq.~(\ref{82}) with $g^{jk}$ gives $ (n-2)^2 \alpha^k\nabla_k R = (n-1)(n-4) \alpha^k \nabla_k G$. Therefore:
\begin{align}
n(n-2)^2 \left (  \frac{\alpha_j \alpha^m}{\alpha^2}-\frac{1}{n}  \delta^m_j  \right )\nabla_m R= (n-1)^2(n-4) \nabla_j G \nonumber
\end{align}
\end{proof}
\end{proposition}
\section{\bf Conformally recurrent Lorentzian manifolds, $\mathbf {n\ge 5}$}
For general Lorentzian manifolds the existence of a covariantly constant symmetric tensor
$h_{ij}$ not proportional to the metric implies reducibility of the holonomy group of the manifold. This result was proven by 
Hall \cite{Hall91} in $n=4$ (see also \cite{ExactSol}), and extended to  $n\ge 5$ by Senovilla (\cite{Senovilla} lemma 3.1), 
Aminova \cite{Aminova} and Galaev \cite{Galaev2}. \\
We recall that the holonomy group is {\em reducible} if it leaves a non trivial subspace of the tangent space invariant; 
it is {\em non-degenerately reducible} (or decomposable) if it leaves a non-degenerate subspace (i.e. such that the restriction of the metric 
is non-degenerate) invariant. These concepts were introduced by Wu for arbitrary signature of the metric \cite{Wu}. Wu showed 
that if the manifold is a connected space-time whose holonomy group is non-degenerately reducible,
then the manifold is locally decomposable.\\
Senovilla's Lemma 3.1 \cite{Senovilla} 
makes the following statement:
Let $(\mathcal M, g)$ be a $n$-dimensional Lorentzian manifold equipped with a covariantly 
constant tensor field $h$ not proportional to the metric. Then the manifold is 
reducible, and it is not non-degenerately reducible only if there 
exists a null covariantly constant vector field which is the unique (up to a 
proportionality constant) constant vector field.\\
Accordingly, there are three possibilities: 1) $\mathcal M$ is decomposable and there are no parallel null
vector fields, 2) $\mathcal M$ is decomposable and there are at least two independent parallel null vector
fields, 3) $\mathcal M$ is not decomposable (and either there are no parallel null vector fields or a
unique one, up to a constant scaling).

\begin{remark}{\quad}\\
In dimension $n=4$, a manifold whose holonomy is type R$_4$ has two independent covariantly constant null 
covector fields $l$ and $n$ with $l^kn_k=1$. Then $z=l+n$ and $t=l-n$ are covariantly constant and $z^k z_k=2$, $t^kt_k=-2$, $z^kt_k=0$.
Then $Z_{ij}=z_iz_j/(z^kz_k)$ and $T_{ij} = t_it_j /|t^kt_k|$ are orthogonal projectors and the metric is decomposable as
$g_{ij}=T_{ij} +Z_{ij}+ (g_{ij}-T_{ij}-Z_{ij})$ i.e. the manifold is 1+1+2 decomposable (see \cite{Hall04} p. 245).
\end{remark}

Aminova  \cite{Aminova} proved that if a locally indecomposable Lorentzian manifold admits a covariantly constant bilinear form 
not proportional to the metric, 
then the manifold admits a covariantly constant null vector field $\beta $ and the space of covariantly constant bilinear forms 
is 2-dimensional and it is spanned by $g$ and by $\beta\otimes\beta $. Locally:
\begin{align} 
h_{ij}=A g_{ij} + B \beta_i\beta_j
\end{align} 
with constants $A$, $B$. The metric takes the form stated by Brinkmann \cite{Brinkmann}:
\begin{align}
ds^2 = -2dx^0\,(Hdx^0 + dx^1 + W_\nu dx^\nu) + G_{\mu\nu} dx^\mu dx^\nu , \quad (\mu,\nu =2,\dots ,n-1), \label{Brink}
\end{align}
where the functions $H$, $W_\nu$ and $G_{\mu\nu}=G_{\nu\mu}$ are arbitrary but independent of $x^1$. These metrics are called {\em Brinkmann waves}. In particular it is a pp-wave if the curvature tensor satisfies the trace condition $R_{ij}{}^{pq}R_{pqlm}=0$ \cite{Lei06a,Lei06b,Lei10}. A Lorentzian manifold is a pp-wave if and only if the metric has the local form Eq.~(\ref{Brink}) with $W_\mu=0$ and $G_{\mu\nu}= \delta_{\mu\nu}$ \cite{Lei06a,Lei06b,Lei10,Sch74};  the function $H$
is usually called the potential of the pp-wave.

\begin{thrm}
A Lorentzian conformally recurrent manifold of dimension $n\ge 5$ with non-zero scalar $C^2$, is decomposable.
\begin{proof}
Suppose that the manifold is non decomposable. 
Then there exists a vector field $\beta_j$ such that $\beta^j\beta_j=0$, $\nabla_i\beta_j=0$. The
tensor $h_{ij}$ in Eq.~(\ref{tensorh}) gains the form
$h_{ij}=Ag_{ij} + B \beta_i\beta_j$. Since $h^k{}_k=1$ it is $A=1/n$. The vectors $\alpha^i$ and $\beta^i$
are eigenvectors of $h_{ij}$ with eigenvalues $\frac{n-3}{2(n-2)}\neq \frac{1}{n}$. Then $\alpha_i\beta^i=0$.
Contraction of $h_{ij}$ with $\alpha^i$ gives: $ \frac{n-3}{2(n-2)} \alpha_j = \frac{1}{n}\alpha_j $.
The factors in the two sides do not match.
\end{proof}
\end{thrm}

It turns out that $C^2=0$ is a necessary condition for a Lorentzian conformally recurrent manifold of dimension $n\ge 5$ to be
non-decomposable.
Lorentzian conformally recurrent manifolds were classified by Galaev \cite{Galaev}, and their metric is that of a pp-wave:
\begin{thrm}[Galaev]
If $\mathcal M$, $n\ge 4$, is  a locally non-decomposable Lorentzian conformally recurrent manifold, 
then either $C_{ijkl}=0$, or $\nabla_m R_{ijkl} = 0$, or locally the metric has the form
\begin{align}
 ds^2 =-2dvdu - (du)^2 \sum_{i=3}^n (x^i)^2 [a(u)+F(u) \lambda_i]  + \sum_{i=3}^n(dx^i)^2 \end{align}
where $a(u)$, $F(u)$ are functions, $\lambda_i $ are real numbers with $\lambda_3 +\dots + \lambda_n=0$. In particular:
$\nabla_i C_{jklm}=0$ if and only if $F(u)$ is constant; $\nabla_i R_{jklm}=0$ if and only if $F(u)$ and $a(u)$ are constants; the Riemann
tensor is recurrent if and only if $F(u)=a(u)$, or $a(u)=0$ or all $\lambda_i=0$; finally, $C_{ijkl}=0$ if and only if $F(u)\lambda_i=0$
for all $i$.
\end{thrm}
In Galaev's theorem the pp-wave potential is $H(\vec x, u)=\sum_{i=3}^n (x^i)^2 [a(u)+F(u) \lambda_i] $. The
Ricci tensor is
\begin{align}
R_{ij} = -\frac{1}{2} \sum_{k=3}^n \frac{\partial^2 H}{\partial x_k^2}\beta_i\beta_j =(n-2) a(u) \beta_i\beta_j
\end{align}
where the vector $\beta_j =\nabla_ju $ is covariantly constant and null. The Ricci tensor is rank-1 and traceless (the scalar curvature being zero) (see \cite{Pod09} sections 2 and 7.1, and \cite{Hall04} page 248). It immediately follows that
\begin{align}
\nabla_kR_{ij} = (n-2)\frac{da}{du}(\nabla_k u)\beta_i\beta_j
=(n-2) a'(u)\beta_k\beta_i\beta_j
\end{align}
From Eq.~(\ref{divC}) for the divergence of the Weyl tensor, we infer that $\nabla_m C_{jkl}{}^m=0$.
\begin{proposition}
A locally non-decomposable Lorentzian conformally recurrent manifold, $n\ge 5$, 
which is neither conformally flat nor locally symmetric, is conformally harmonic.
\end{proposition}

\section{\bf Algebraic classification of Lorentzian conformally recurrent manifolds}
If $\alpha^k\alpha_k<0$, the tensor $E_{ij}$ defined in Eq.~(\ref{eq2.15a}) is referred to as the ``electric tensor''
associated to the Weyl tensor. In general Lorentzian $n=4$ manifolds, the Weyl tensor is completely described by the electric tensor and
a ``magnetic tensor'' \cite{BertschHam, Mat53,Mcl94}. The decomposition was extended to $n\ge 5$ by Senovilla \cite{Sen01}, Hervik et al. \cite{HeOrWy2013} and Ortaggio et al. \cite{Ortaggio13}:
given a time-like vector $u_k$ ($u^2=-1$) introduce the tensors
$\theta_{kl} = g_{kl} +u_ku_l$ and $E_{kl} =u^ju^mC_{jklm}$. The Weyl tensor is the sum
$C=C_+ + C_-$, with electric and magnetic components:
\begin{align}
&(C_+)^{jk}{}_{ml} = \theta^k{}_s \theta^s{}_l E^j{}_m +4u^{[j} u_{[l} E^{k]}{}_{m]}\\
&(C_- )^{jk}{}_{ml} = 2 \theta^{jr} \theta^{ks} C_{rsp[l} u_{m]} u^p+ 2 u_p u^{[j} C^{k]prs}\theta_{rl} \theta_{sm}.
%\label{eq2.18}
\end{align}
The Weyl tensor is purely electric ($C_-=0$) if and only if
$u_iu_m C_{jkl}{}^m + u_ju_m C_{kil}{}^m + u_k u_m C_{ijl}{}^m = 0 $ (theorem 3.5 in ref.\,\cite{HeOrWy2013}).
Because of Eq.~(\ref{eq2.17}) we may assert:

\begin{proposition}\label{Theorem 3.5}
On a Lorentzian conformally recurrent manifold $\mathcal M$ 
of dimension $n\ge 5$, with time-like recurrence vector ($\alpha^2 <0 $) the Weyl tensor is purely electric.
\end{proposition}

In $n=4$ the algebraic classification of the Weyl tensor on Lorentzian conformally recurrent manifolds
was obtained by Mc Lenaghan and Leroy \cite{LenaghLer}. On $n\ge 5$ we note that
Eq.~(\ref{eq2.17}) matches with type $II_d$ of table 1 of \cite{Ortaggio09}, that classifies
the Weyl tensors on $n\ge 5$ Lorentzian manifolds, generalising the Bel-Debever scheme
based on null vectors. Then:

\begin{proposition}\label{Theorem 3.2}
On a Lorentzian conformally recurrent manifold $\mathcal M$, $n\ge 5$, with null recurrence vector, $\alpha^2 = 0$,
the Weyl tensor is algebraically special at least as type II$_d$
of the Bel-Debever-Ortaggio classification.
\end{proposition}
We can be more specific (see \cite{Ortaggio13}). For type$-N$ spaces it is $C_{pklm}C_q{}^{klm}=0$, then $C^2=0$ (pp-waves are in this category). \\
In type$-III$ spaces
$C_{pklm}C_q{}^{klm}=\overline \psi l_pl_q$, where $l$ is the standard null vector of a Lorentzian basis; again it is $C^2=0$.
Recurrency of the Weyl tensor implies that $\nabla_p l_q =  (\alpha_p -\frac{1}{2}\nabla_p \log \overline\psi)l_q$. Therefore, a type$-III$
conformally recurrent manifold admits a recurrent null vector field. The metric, in Walker's
coordinates \cite{Walker49}, belongs to the Kundt's class \cite{Pod09}:
$$ds^2 =2dudv  +a_i(\vec x, u)dx^idu+H(v,\vec x,u)du^2+g_{ij}(\vec x, u) dx^idx^j.$$
If $a_i$ is locally closed, the vector $l_p$ can be rescaled to a null constant vector, and the metric becomes of
Brinkmann's type.\\
For types equal or higher than $III$ it is always $C^2=0$, as guaranteed by a theorem by Hervik \cite{Hervik2011}.
\section{\bf Examples of conformally recurrent manifolds}
The first example of a conformally recurrent manifold was given by Roter in \cite{Roter0} and further discussed in \cite{Roter1,Roter2,Roter3},
$ds^2 = Q\, (dx^1)^2 + k_{\mu\nu} dx^\mu dx^\nu + 2 dx^1dx^n$,  $1<\mu,\nu <n$,
$Q=[A(x^1)p_{\mu\nu} + k_{\mu\nu}]x^\mu x^\nu $ with tensors $p$ and $k$ subject to restrictions.
The manifold is both conformally recurrent and Ricci-recurrent, possibly with different recurrence vectors.
The same metric structure was investigated and classified according to choices of the weights by Grycak and Hotlo\'s \cite{Grycak82}.
Derdzi\'nski \cite{Derdzinski} studied the metric $g_{11}=-2\epsilon$, $g_{ij}=\exp [F_i (x^1,x^2)]$ if $i+j=n+1$, and $g_{ij}=0$
otherwise, with certain functions $F_i$, and periodicity conditions.

Here we give another example. Consider the Brinkmann-type metric:
\begin{align}
 ds^2=p(x^1)q(x^3) (dx^1)^2 + 2 dx^1dx^2 + (dx^3)^2 +\dots +(dx^n)^2. \label{esempiobello}
 \end{align}
The non-zero Christoffel symbols (up to symmetries) are:
$ \Gamma_{11}^2 = \tfrac{1}{2}p'(x^1)q(x^3) $, $\Gamma_{13}^2 =-\Gamma_{11}^3 =\tfrac{1}{2} pq' $.
It follows that the non-zero components of the Riemann tensor are those related to
$ R_{1313} =\frac{1}{2}pq''$, and the only non-zero
component of the Ricci tensor is $R_{11}=\frac{1}{2}pq''$. The curvature scalar is zero.\\
It is easy to check that the manifold is Ricci-recurrent, $\nabla_i R_{jk} = \alpha_i R_{jk} $, with recurrence vector
$$\alpha_i = \left ( \frac{p'(x^1)}{p(x^1)}, 0 , \frac{q'''(x^3)}{q''(x^3)}, 0,\dots ,0 \right ).$$
Further, the manifold is recurrent with the same recurrence vector: $\nabla_i R_{jklm}=\alpha_i R_{jklm}$.
Therefore the Weyl tensor is recurrent, i.e. the manifold is a conformally recurrent manifold.
Since $\alpha $ is a gradient, by \eqref{eq3} the metric is Weyl semi-symmetric (for the same reason
it is Ricci semi-symmetric and semi-symmetric).\\
The form of the metric gives $\alpha^k\alpha_k = (q'''/q'')^2\ge 0$.\\
The simple forms of the Riemann and Ricci tensors in the defining frame Eq.~(\ref{esempiobello}),
imply the tensorial identity $\alpha^k\alpha^m R_{jklm} = \alpha^i\alpha_i R_{jl}$.
Since the curvature scalar is zero, it follows that $\alpha^k\alpha^m R_{km}=0$.

If $\alpha^k\alpha_k=0$ then $\alpha = (p'/p, 0,0,\dots)$ and $R_{ij}\propto \alpha_i\alpha_j$ in the defining
frame. Then in any frame the Ricci tensor is rank-one: $R_{ij}= \frac{1}{2}\alpha_i\alpha_j\,
(pq'')(p'/p)^{-2}$. This has the following consequences:
$\alpha^m R_{jklm}= \nabla^m R_{jklm} = \nabla_k R_{jl}-
\nabla_j R_{kl} = \alpha_k R_{jl}-\alpha_j R_{kl}=0$ and, after simple calculations,
$\alpha^m C_{jklm}=0$. For a Lorentzian metric ($|pq|>2$), the last result characterizes the
manifold as type $II'_{abd}$ in the Bel-Debever-Ortaggio classification \cite{Ortaggio09}.
\section*{Acknowledgements}
We dedicate this paper to the memory of  Dr. Wies{\l}aw Grycak, whose research on conformally recurrent manifolds much inspired ours.\\
We cordially thank the Referees for their accurate and useful suggestions.

\end{document}